\documentstyle[12pt]{article}

\newcommand{\be}{\begin{equation}}
\newcommand{\ee}{\end{equation}}
\newcommand{\re}[1]{(\ref{#1})}
\newcommand{\hslq}{U_q\widehat{sl}_2}
\newcommand{\slq}{U_qsl_2}
\newcommand{\slh}{U_hsl_2}
\newcommand{\sle}{U_\varepsilon sl_2}
\newcommand{\hslh}{U_h\widehat{sl}_2}
\newcommand{\hsle}{U_\varepsilon\widehat{sl}_2}
\newcommand{\ot}{\otimes}
\newcommand{\mod}[1]{|#1|}
\newcommand{\half}[1]{{#1 \over 2}}
\newcommand{\binom}[2]{\left[\begin{array}{c} #1 \\ #2 \end{array}\right]}
\newcommand{\q}{{q^{-2}}}

\newcommand{\Ad}{{\rm Ad}}
\newcommand{\ad}{{\rm ad}}
\newcommand{\Li}{{\rm Li_2}}
\newcommand{\ab}[1]{\rm{\bf{#1}}}

\newcommand{\ve}{\varepsilon}
\newcommand{\lo}{{\lambda_1}}
\newcommand{\lt}{{\lambda_2}}

\begin{document}

\rightline{SPhT-94-146}

\begin{center}

{\large {\bf On the universal $R$-matrix of  $U_q\widehat{sl}_2$

	at roots of unity}}
	\footnote{This work is supported by NATO linkage grant
	LG 9303057}

\vspace{24pt}

		T.Hakobyan
	\footnote{Permanent address: {\sl Yerevan Physics Institute,
	375036 Yerevan, Armenia,}

	 e-mail: {\sl hakob@vxc.yerphi.am
	or hakob@uniphi.yerphi.am}}
and	A.Sedrakyan
	\footnote{Permanent address: {\sl Yerevan Physics Institute,
	375036 Yerevan, Armenia,}

	{\sl e-mail: sedrak@vxc.yerphi.am
	or sedrak@uniphi.yerphi.am}}

\vspace{10pt}

{\sl Service de Physique Theorique de Saclay,} \\
{\sl F-91191, Gif-sur-Yvette, France.}

\end{center}

\vspace{2cm}

\begin{abstract}
We show that the action of universal $R$-matrix of  affine $U_qsl_2$ quantum
algebra, when $q$ is a root of unity, can be renormalized by some scalar
factor to give a well defined nonsingular expression, satisfying
Yang-Baxter equation. It reduced to intertwining operators
of all representations, corresponding to Chiral Potts, if the
parameters of these representations lie on well known algebraic curve.

We also show that affine $U_qsl_2$ for $q$ is a root of unity form
the autoquasitriangular Hopf algebra
in the sence of Reshetikhin.
\end{abstract}

\newpage

\section{Introduction}

The intertwining operators of quantum groups (\cite{Dr1,J1,J2,J3}) lead to
solutions of
Yang-Baxter equation, which play the crucial role in
two dimensional field theory and integrable statistical systems
(\cite{J3,PS}). It is well known that the most of them can be obtained
from the universal $R$-matrix (\cite{Dr1}) for a given quantum group:
the solutions of spectral parameter dependent Yang-Baxter equation
can be obtained from the universal $R$-matrix of affine quantum groups
(\cite{Tol2}) and the solutions of non-spectral parameter dependant Yang-Baxter
equations can be obtained from the universal $R$-matrix  of finite quantum
groups.

The situation is not the same for the case, when the parameter $q$
of quantum group is a root of unity.

In this case the center of
quantum group is larger and new type of representations appear, which have
no a classical analog (\cite{PS,ACh1,ACh2,Kac}).
It was shown in \cite{BS,BKMS} that the cyclic representations
lead to solutions of Yang-Baxter equation with a spectral parameter,
lying on some algebraic curve. These solutions correspond to Chiral Potts
Model
(\cite{YMcPTY,McPTS,BPY}) and its generalizations (for quantum groups
$U_q\widehat{sl}_n$).

The formal expression of the universal $R$-matrix fails in this
case: it have a singularities at $q$ is a root of unity. Recently
in \cite{Resh2}
Reshetikhin introduced the notion of autoquasitriangular Hopf algebra
to avoid these singularities. He treated the $\slq$ case.

The main goal of this paper is to show that after suitable renormalization
by scalar factor the universal $R$-matrix produces $R$-matrices for
concrete representations.

In section 2, we consider the universal $R$-matrix on Verma modules
of $\slq$ for $q$ is a root of unity. We prove that it is well defined
and make a connection with the $R$-matrix of autoquasitriangular Hopf
algebra, founded by Reshetikhin.

In section 3 we consider the algebra $\hslq$ at roots of unity.
We found the central elements of its Poincar\'e-Birkoff-Witt (PBW) basis,
generalizing the results of
\cite{Kac} for affine case. It appears that new type of central elements
appear for some imaginary roots, which have no analog for finite quantum
groups.
After this we prove the autoquasitriangularity of $\hslq$, generalizing
the results of \cite{Resh2} for affine case. Then we consider the action of
affine universal
$R$-matrix on $\hslq$- and $\slq$-Verma modules. On $\hslq$-Verma modules
it is well defined. For $\slq$-Verma modules (evaluation representation)
we renormalize its expression
by scalar factor to exclude the singularities. The rest part lead to solutions
of infinite dimensional spectral parameter dependent Yang-Baxter equation.
We showed that under the certain condition this $R$-matrix can be restricted
to semicyclic representations, giving the Boltsmann weights of Chiral Potts
model, corresponding to such type representations, which was considered in
\cite{GSR-A,GS,IU,HS1}. The condition, mentioned above, is on the parameters of
representations: they must lie on well known algebraic curve. It is
integrability condition of Chiral Potts model.

In the last section we made same type suggestion for the cyclic
representations.

\section{The $\slq$ case}
\subsection{The Universal $R$-matrix on Verma Modules at root of unity}

The quantum group $\slq$ is a ${[}q,q^{-1}{]}$-algebra, generated by the
elements $E,\ F,\ K$ with the following relations between them

\begin{eqnarray} \label{def_sl}
&{[}K,K^{-1}{]}=0 &
{[}E,F{]}= \frac{K-K^{-1}}{q-q^{-1}} \\
&KEK^{-1}=q^2E & KFK^{-1}=q^{-2}F, \nonumber
\end{eqnarray}

On $\slq$ there is a Hopf algebra structure with  comultiplication
$\Delta: \ \slq\rightarrow\slq\ot\slq$ defined by
\begin{eqnarray*}
\label{comul_sl}
\Delta(K)=K\ot K & \Delta(E)=E\ot 1 +K^{-1}\ot E &
\Delta(F)=F\ot K+1\ot F
\end{eqnarray*}

We denote $K=q^H, \ q=e^\hbar$, as usually, and consider the $[[h]]$-algebra
$\slh$ with the same defining relations. $\slh$ is a quasitriangular Hopf
algebra, i.e. it possess
the universal $R$-matrix $R\in\slh\ot \slh$ connecting the comultiplication
$\Delta $ with the opposite comultiplication $\Delta^\prime=\sigma\circ\Delta$,
where $\sigma(x\ot y):=y\ot x$:

\begin{equation}
\label{inter}
\Delta^\prime(a)=R\Delta(a)R^{-1}, \qquad \forall a\in \slq
\end{equation}

It satisfies the quasitrialgularity relations
\begin{eqnarray}
\label{quasi}
(\Delta\ot 1)R=R_{13}R_{23} & (1\ot \Delta)R=R_{13}R_{12}
\end{eqnarray}
and Yang-Baxter equation (\cite{Dr1})
\begin{equation}
\label{YB}
	R_{12}R_{13}R_{23}=R_{23}R_{13}R_{12}
\end{equation}
Here we used the usual notation: if $R=\sum_i a_i\ot b_i$, $a_i,b_i\in \slq$,
then
$$
R_{12}=\sum_i a_i\ot b_i\ot 1 \ \ R_{13}=\sum_i a_i\ot 1\ot b_i \
\ R_{23}=\sum_i 1 \ot a_i \ot b_i
$$
The explicit expression of $R$ in terms of formal power series is
\begin{eqnarray}
\label{sl_rmat}
R=\exp_{q^{-2}}({(q-q^{-1})(E\ot F)})q^{{1\over 2}H\ot H}
\end{eqnarray}
where the $q$-exponent is defined by
$\exp_q(z)=\sum_{n\ge 0}{z^n\over(z)_q!}$,
$(z)_q:=\frac{1-q^n}{1-q}$.

Note, that to be precise, $\slq$ is not a quasitriangular Hopf algebra,
because the term $q^{H\ot H}$ in \re{sl_rmat} do not belong to
$\slq\ot\slq$, but it
is an autoquasitriangular Hopf algebra (\cite{Resh2}). The latter
is a Hopf algebra $A$, where the condition \re{inter} is generalised by
$$
\Delta^\prime=\hat{R}(\Delta),
$$
where $\hat{R}$ is an automorphism  of $A\ot A$ (not inner, in general).
So, although \re{sl_rmat} is ill defined on $\slq$, but the
action
\begin{equation}
\label{hR}
\hat{R}(a)=RaR^{-1},
\end{equation}
where $a\in \slq\ot\slq$ is till well defined.

For two representations of $\slq$ $V_1$ and $V_2$ one can consider two
$\slq$-actions on $V_1\ot V_2$ by means of both comultiplications $\Delta$
and $\Delta^\prime$. If $R$ is defined on $V_1\ot V_2$, then both
$\Delta$- and $\Delta^\prime$-actions are equivalent via intertwining
operator $R_{V_1\ot V_2}=R|_{V_1\ot V_2}$. For general $q$ the restriction
of \re{sl_rmat} on tensor product of two irreducible representations
(in general, of any highest weight representations) is well defined.
And all solutions of Yang-Baxter equation \re{YB}, having $\slq$-
symmetry in sense of \re{inter} can  be obtained from the universal
$R$-matrix \re{sl_rmat} in such  way.

The situation is different for $q$ being a root of unity.
In this case the singularities appear in the formal expression of $R$.

Recall that for
$q=\exp\left(\frac{2\pi i}{N^\prime}\right)$
the elements $F^N,\ F^N,\ K^N$,
where $N=N^\prime$ for odd $N^\prime$ and $N={N^\prime\over 2}$ for even
$N^\prime$, belong to the center of $\slq$. In irreducible representations
they are multiples of identity.
Recall that every $N$-dimensional irreducible representation is
characterized by the values $x,y,z$
of these central elements (and also by the value of $q$-deformed Casimir
$c=\frac{Kq+K^{-1}q^{-1}}{q-q^{-1}}+EF$, which for the fixed $x,y,z$ can have
in general $N$ discrete values (\cite{Kac})).

Although the expression of $R$-matrix \re{sl_rmat} of  $\slq$ has
singularities for $q^{N^\prime}\to 1$ in all terms ${1\over (n)_\q!}E^n\ot F^n$
for $n\ge N$,
its restriction on tensor product of Verma modules $M_\lo\ot
M_\lt$ is well defined.

Recall that $M_\lambda$ is formed by the
basic vectors $v_m^{\lambda}$, $m=0,1,\dots$, satisfying

\begin{eqnarray*}
\label{Verma}
 Ev_0^\lambda=0 & Fv_m=v_{m+1}^\lambda & Hv_0=\lambda v_0, \quad\lambda\in
\ab{C}
\end{eqnarray*}

To consider the action of $R$ on $M_\lo\ot M_\lt$
we use the formula, which can be obtained from the defining relations
\re{def_sl} (\cite{Kac}):

$$
[E^n,F^s]=\sum_{j=1}^{\min(n,s)}\binom{n}{j}\binom{s}{j}[j]!F^{s-j}
        \left( \prod_{r=1}^{j}[H+j-n-s+r] \right) E^{n-j},
$$
where
$$
\binom{a}{b}=\frac{[a]!}{[b]![a-b]!} \ \ {\rm and} \ \
{[}n{]}=[n]_q=\frac{q^n-q^{-n}}{q-q^{-1}}.
$$

So, for $n>s$  ${E^n \over (n)_\q !} v_s^\lambda=0$ and for  $n\le s$:

$$
  {E^n \over (n)_\q !}v_s^\lambda=q^\half{n(n-1)}\binom{s}{n}\prod_{r=1}^n
        [\lambda-s+r] v_{s-n}^\lambda
$$

The $q$-binomial $\binom{s}{n}$ has a non-infinity limit
for $q^{N^\prime}\to 1$. So,

\begin{eqnarray}
\label{r_verma1}
R(v_s^\lo\ot  v_{s^\prime}^\lt ) &=&
        \sum_{n=0}^s q^\half{(\lambda_1-2s)(\lambda_2-2s^{'})}q^\half{n(n-1)}
	 (q-q^{-1})^n \binom{s}{n}
	\nonumber\\
        &\times&\prod_{r=1}^n [\lambda_1-s+r]
        v_{s-n}^\lo\ot v_{s^\prime+n}^\lt
\end{eqnarray}

is well  defined for $q$ is a root of unity.

\subsection{The connection with Reshetikhins $R$-matrix of
autoquasitriangular Hopf algebra}

This $R$-matrix can be presented in another form by using  resent
results of  Reshetikhin (\cite{Resh2}). He used an  asymptotic formula
for $q$-exponent in the limit $q^{N^{'}}\to 1$
to bring out multiplicatively  singularities from
$\exp_q((q-q^{-1})E\ot F)$. The expression of universal $R$-matrix in this
limit then acquires the form:

\begin{eqnarray*}
\label{rm_1}
R&=&\exp\left(\frac{1}{2N^2\hbar}\Li(E^N\ot F^N)\right)(1-E^N\ot
F^N)^{-\half{1}}
	\nonumber\\
 &\times&       \prod_{m=0}^{N-1}(1-\ve^mE\ot F)^{-\frac{m}{N}}
        q^{\half{1}H\ot  H} \cdot O(\hbar)
\end{eqnarray*}

Here $q=\exp(\hbar) \ve,  \ \ve=\exp(\frac{2\pi i}{N^\prime})$ and
$\Li(x)=-\int\limits_0^x\frac{\ln(1-y)}{y}dy$ is a dilogaritmic function.

Recall that although the elements
\begin{eqnarray}
\frac{E^N}{(N)_\q!}, \quad \frac{F^N}{(N)_\q!} \quad  {\rm and} \  H
\nonumber
\end{eqnarray}
don't belong to
$\slq$ for $q^{N^\prime}\to 1$, but their adjoint actions
$$
\ad(x)a={[}x,a{]},\ \Ad(\exp(x))a=\exp(x)a\exp(-x)=\exp(\ad(x))a
$$
on $\slq$ are
well defined in this limit and give rise to some derivations (\cite{Kac}).
Let's denote them by $e,\ f$ and $ h$ correspondingly.

The element $\frac{1}{2\hbar N^2}\Li(E^N\ot  F^N)$ in  the  exponent
of \re{rm_1} in the adjoint representation also acts on $\slq\ot\slq$
as a derivation in the  limit $\hbar\to 0$ .
It can be expressed by means of the derivations
$e$ and  $f$ as follows:
\begin{eqnarray*}
\label{ENFN1}
\lim_{\hbar\to 0}\ad\left(\frac{1}{2\hbar N^2}\Li(E^N\ot F^N)\right)&=&
c_{N^\prime}\frac{\ln(1-E^N\ot F^N)}{E^N\ot F^N}\nonumber\\
	&\times& (e\ot F^N + E^N\ot f),
\end{eqnarray*}
where
\begin{eqnarray}
\label{cN}
c_{N^\prime}= \left\{ \begin{array}{lll}
		-(1-\ve^{-2})^{-N} & {\rm for \ odd} \ N & (N^\prime=N)\\
		(-1)^N(1-\ve^{-2})^{-N} & {\rm for \ even} \ N & (N^\prime=2N)
		\end{array}
		\right.
\end{eqnarray}
Note, that
$$
\Ad(\ve^{\half{1}H\ot H})=\ve^{\half{1}(h\ot H +H\ot h)}=
        K^{1\ot\half{1}h}\ot K^{\half{1}h\ot 1}
$$
is well defined in the adjoint representation.

So, one can write down the automorphism $\hat{R}$  \re{hR} in the limit
$\hbar\to  0$, obtained in \cite{Resh2}, in the following form
\footnote{For quantum groups one  can   introduce 4 equivalent
	co\-mul\-ti\-pli\-ca\-tions:  $\Delta_q$, $\Delta_q^\prime$,
	$\Delta_{q^{-1}}$,
	$\Delta_{q^{-1}}^\prime$ \cite{Tol2}. In \cite{Resh2} the
	comultiplication
	$\Delta_{q^{-1}}^\prime$ had been used as a basic one. So, $R$-matrix,
	used there, is $R_{q^{-1}}^{-1}$  in our notations
	and  differs from $\Delta_q$-case
	used here  by permutation  of $q$-exponent and $q^{\half{1}H\ot H}$.
	 }

\begin{eqnarray}
\label{rm_2}
\hat{R}&=& \prod_{m=0}^{N-1}\Ad\left((1-\ve^mE\ot F)^{-\frac{m}{N}}\right)
	\nonumber\\
	&\times&\exp\left(-(1-\ve^{-2})^{-N}\frac{\ln(1-E^N\ot F^N)}{E^N\ot F^N}
	(e\ot F^N+ E^N\ot f)\right)\\
	&\times& K^{1\ot\half{1}h}\ot K^{\half{1}h\ot 1}
\nonumber
\end{eqnarray}

Let us now consider the restriction of \re{rm_2} on the quotient
algebra obtained
from $\slq$ by factorisation on the  ideal, generated by $E^N$, i.e.
impose $E^N=0$.  Although this ideal is not stable with respect to
derivations $e,f,h$, it is easy to see that it is stable with respect
to $\hat{R}$.

Moreover,  the left  $\slq\ot\slq$-module
$$
I_{\lambda_1,\lambda_2}=(\slq\ot I_\lt)\bigoplus (I_\lo\ot\slq)
$$
is also  stable  with  respect to $\hat{R}$. Here we denoted by  $I_\lambda$
the  left $\slq$-module, generated  by $E$ and  $(K-\ve^\lambda)$. This fact
allows  to restrict \re{rm_2}  on  Verma  modules, because we have the left
$\slq$-module equivalence
$$
(\slq\ot\slq) / I_{\lo,lt}\cong M_{\lo}\ot M_\lt
$$

So,  one can derive   from \re{rm_2} the restriction of  $R$ on  this
factormodule is given by  the  multiplication on
\footnote{Note that both $h$ and $e$ are well defined on $M_\lambda$
in contrast to $f$}

\begin{eqnarray}
\label{r_verma2}
R&=& \prod_{m=0}^{N-1}\left((1-\ve^mE\ot F)^{-\frac{m}{N}}\right)
\nonumber\\
       &\times& \exp\left((1-\ve^{-2})^{-N}
       (e\ot F^N)\right) K^{1\ot\half{1}h}\ot K^{\half{1}h\ot 1}
\end{eqnarray}
This expression is another form of expression of universal
$R$-matrix \re{sl_rmat} on Verma modules and coincide with \re{r_verma1}.

\section{The case of affine $\hslq$}
\subsection{The PBW basis and the universal $R$-matrix}
 The affine quantum universal enveloping algebra $\hslq$ is a
$[q,q^{-1}]$-Hopf algebra, generated by elements
$E_i:=E_{\alpha_i}$, $F_i=F_{\alpha_i}$, $K_i=q^{H_i}$, $i=0,1$ and  $q^d$
with defining relations (\cite{J1}):

\begin{eqnarray} \label{def_hsl}
&[q^{H_i},q^{H_j}{]}=0 \qquad q^dq^{H_i}=q^{H_i}q^d \qquad
{[}E_i,F_j{]}=\delta_{ij} {[}H_i{]}_q
\nonumber\\
&q^{H_i}E_jq^{-H_i}=q^{a_{ij}}E_j \qquad q^{H_i}F_jq^{-H_i}=q^{-a_{ij}}
\qquad q^dE_1q^{-d}=qE_1
\\
& q^dF_1q^{-d}=q^{-1}E_1 \qquad q^dE_0q^{-d}=E_0 \qquad q^dF_0q^{-d}=F_0
\nonumber\\
& (\ad_qE_i)^{1-a_{ij}}E_j=0 \qquad (\ad_qF_i)^{1-a_{ij}}F_j=0
 \nonumber
\end{eqnarray}
and comultiplication
\begin{eqnarray}
\label{comul}
&\Delta(q^{H_i})=q^{H_i}\ot q^{H_i}
&\Delta(q^{d})=q^{d}\ot q^{d}  \nonumber\\
&\Delta(E_i)=E_i\ot 1 +q^{-H_i}\ot E_i
&\Delta(F_i)=F_i\ot q^{H_i}+1\ot F_i \nonumber
\end{eqnarray}

Here we use the $q$-deformed adjoint action
$(\ad_qx)y:=\sum_i x_iys(x^i)$, where $\Delta(x)=\sum_ix_i\ot x^i$ and
$s: \ \hslq\to\hslq$ is antipode of $\hslq$, defined by
\begin{eqnarray*}
s(E_i)=-K_iE_i & s(F_i)=-F_iK_i^{-1} & s(K_i)=K_i^{-1}
\end{eqnarray*}
Also we denoted by $a_{ij}$ the Cartan matrix of affine $\hat{sl}(2)$ Lie
algebra
$$a_{ij}=\left( \begin{array}{rr} 2 & -2
				\\ -2 & 2 \end{array} \right)$$
Let's denote by $c$ the central element $c=H_1+H_2$.

Define on  $\hslq$ an antiinvolution $\iota$ by
\begin{eqnarray*}
\label{antiinv}
\iota(K_i)=K_i^{-1} & \iota(E_i)=F_i &\iota(F_i)=E_i \quad\iota(q)=q^{-1}
\end{eqnarray*}

As above, denote by $\hslh$ the $[[h]]$-algebra with the same relations but
the elements $H_i$ instead of $q^{H_i}$.

The PBW basis of $\hslh$ is formed by  elements $H_i,\ d$,
$E_{\alpha_i+n\delta}$, $F_{\alpha_i+n\delta}$, $E^\prime_{n\delta}$ and
$F^\prime_{n\delta}$, which are inductively defined by the relations

\begin{eqnarray}
\label{PBW}
&E_{\alpha_0+n\delta}=(-1)^n(\ad_{E^\prime_\delta})^nE_0  \qquad
E_{\alpha_1+n\delta}=(\ad_{E^\prime_\delta})^nE_1  \nonumber\\
& E^\prime_{n\delta}=[2]^{-1}(E_{\alpha_0+(n-1)\delta} E_1 -
	q^{-2} E_1 E_{\alpha_0+(n-1)\delta})   \\
& F_{\alpha_i+n\delta}=\iota(E_{\alpha_i+n\delta}) \qquad
F^\prime_{n\delta}=\iota(E^\prime_{n\delta}) \nonumber
\end{eqnarray}

The expression of the universal $R$-matrix of $\hslh$ is simpler if one
redefine $E^\prime_{n\delta}$ and $F^\prime_{n\delta}$ by means of Schur
polynomials (\cite{Tol2}):

\begin{eqnarray*}
\label{Shur}
&E^\prime_{n\delta}=\sum\limits_{\stackrel{ 0<k_1<\dots<k_m}
	{k_1p_1+\dots+k_m p_m=n }}
        \frac{(q^2-q^{-2})^{\sum p_i-1}}{p_1!\ldots p_m!}
        (E_{k_1\delta})^{p_1}\ldots (E_{k_m\delta})^{p_m}\\
&F^\prime_{n\delta}=\iota(E^\prime_{n\delta}) \nonumber
\end{eqnarray*}

In order to rewrite all the relations between \re{PBW} in compact form
it is suitable to change slightly the basis as follows:

\begin{eqnarray}
\label{compar}
E_{\alpha_0+n\delta}=(-1)^nq^{-2n}x_{n+1}^-k^{-1} &
E_{\alpha_1+n\delta}=(-1)^nq^{-(c+2)n}x_n^+ \\
E^\prime_{n\delta}=(-1)^n\frac{q^{-\half{c}n-2n}}{q^2-q^{-2}}\psi_nk^{-1} &
E_{n\delta}=(-1)^n\frac{q^{-(\half{c}+2)n}}{[2]}a_n
\nonumber
\end{eqnarray}

Then the elements $x^\pm_n, (n\in \ab{Z})$, $a_k,  (k\in \ab{Z},\ k\ne0)$,
$\psi_m, \varphi_{-m},\ (m\ge 1)$ and $\psi_0=\varphi_0^{-1}=k$ satisfy
the following relations

\begin{eqnarray}
\label{Drin}
&{[}a_m,a_n{]} = \delta_{m,-n}\frac{{[}2m{]}{[}mc{]}}{m} \qquad {[}a_m,k{]}=0
\nonumber\\
&kx_m^\pm k^{-1}=q^{\pm 2}x_m^\pm \qquad {[}a_m,x_n^\pm{]}=
	\pm \frac{{[}2m{]}}{m}q^{\mp{\mod{m}\over 2}c} x^\pm_{m+n} \nonumber\\
&x_{m+1}^\pm x_n^\pm-q^{\pm 2}x_n^\pm x_{m+1}^\pm =
	q^{\pm 2}x_m^\pm x_{n+1}^\pm-x_{n+1}^\pm x_m^\pm \\
&{[}x_m^+,x_n^-{]}=\frac{1}{q-q^{-1}}(q^{\frac{c}{2}(m-n)}\psi_{m+n}-
	q^{-\frac{c}{2}(m-n)}\varphi_{m+n}) \nonumber\\
&\sum_{m=0}^\infty \psi_mz^{-m}=k\exp\left({(q-q^{-1})\sum_{m=1}^\infty
a_m z^{-m}}\right)
\nonumber\\
&\sum_{m=0}^\infty \varphi_{-m}z^{m}=k^{-1}\exp\left({-(q-q^{-1})
\sum_{m=1}^\infty a_{-m} z^{m}}\right)
\nonumber
\end{eqnarray}
These relations had been introduced by Drinfeld in \cite{Dr2}
and define another realization of affine algebra $\hslq$.
The antiinvolution $\iota$ in this notations is
\begin{eqnarray*}
\label{antiinvx}
\iota(x^\pm_n)=x^\mp_{-n} & \iota(\psi_n)=\varphi_{-n} &
\iota(a_n)=a_{-n} \qquad \iota (q)=q^{-1}
\end{eqnarray*}

We choose the normal ordering of positive root system $\Delta_+$ of $\hslq$
as follows:
\begin{equation}
\label{normal}
\alpha_0, \alpha_0+\delta,\ldots, \alpha_0+n\delta,\ldots
	\delta,2\delta,\ldots,n\delta,\ldots,
	\alpha_1+n\delta,\ldots, \alpha_1+\delta,\ldots,\alpha_1
\end{equation}
Then the universal $R$-matrix has the form \cite{Tol2}:

\begin{eqnarray}
\label{rmatrix}
R&=&\left(\prod_{n\ge 0}\exp_\q((q-q^{-1})(E_{\alpha_0+n\delta}
        \ot F_{\alpha_0+n\delta}))\right)
\nonumber\\
&\times&        \exp\left(\sum_{n>0}\frac{n E_{n\delta}\ot F_{n\delta}}
	{q^{2n}-q^{-2n}}\right)
\\
&\times& \left(\prod_{n\ge 0}\exp_\q((q-q^{-1})(E_{\alpha_1+n\delta}
        \ot F_{\alpha_1+n\delta}))\right)
        q^{{1\over 2}H_0\ot H_0+c\ot d+d\ot c}
\nonumber,
\end{eqnarray}
where the product is given according to the normal order \re{normal}.

\subsection{ $\hslq$ at roots of unity}

For $q$ being a root of unity ($q=\ve$, $\ve=
e^{\frac{2\pi i}{N^\prime}}$) the center of $\hslq$ is enlarged by the
$N$-th power of the root vectors, as for finite quantum groups:

\begin{eqnarray}
\label{center}
[E_\gamma^N,x]=0 & [F_\gamma^N,x]=0 & [K_i^N,x]=0,
\end{eqnarray}
where $\gamma\in \Delta_+:=\{\alpha_i+n\delta,m\delta|n\ge0,m>0\}$ and
$x\in\hsle$.

These conditions for the simple roots $\gamma=\alpha_i$ can be proven
by using the defining relations of Cartan-Weyl basis \re{def_hsl} as
for finite quantum algebras it had been done in \cite{Kac}.
Indeed, using
$$
\Delta(E_i^N)=K_i^{-N}\ot E_i^N+E_i^N\ot 1,
$$
recalling that $q$-deformed adjoint action $\ad_q$ is a
$\hslq$-representation:
$$
\ad_q(ab)c=\ad_q(a)ad_q(b)c, \ \forall a,b,c\in\hslq
$$
and using Serre relations in \re{comul}, we obtain for $i\ne j,\ N\ge3$:
$$
[E_i^N,E_j]=\ad_q(E_i^N)E_j=(\ad_q(E_i))^NE_j=0
$$
Other commutations in \re{center} for $\gamma=\alpha_i$ can be
verified easily.

To carry out \re{center} for other roots one can try to use the isomorphism,
induced by the $q$-deformed Weyl group. In affine case it had considered
in \cite{Lev}. But it is easier to use the
symmetries of Drinfeld realization of $\hslq$ directly. It is easy
to see from \re{Drin} that  the operation $\omega_\pm$ on $\hslq$ defined by

\begin{eqnarray}
\label{sigma}
 \omega_\pm(x_m^\pm)=x_{m\pm1}^\pm & \omega_\pm(a_m)=a_m & \omega_\pm(q)=q
\nonumber\\
 \omega_\pm(\psi_n)=q^c\psi_n & \omega_\pm(\varphi_n)=q^{-c}\varphi_n &
 \omega_\pm(c)=c
\end{eqnarray}
is an algebra automorphism. As the roots can be obtained by
applying $\omega_\pm$ from the simple ones, we finished the proof.

In addition to this, the elements $E_{kN\delta}$, $F_{kN\delta}$ are
central for $k\in{\bf N}_+$. This can be seen from \re{Drin} and
\re{compar}. These central elements have no analog for finite
algebras.

The adjoint action of $\frac{E_\gamma^N}{(N)_\q!}$,
$\frac{F_\gamma^N}{(N)_\q!}$, $\gamma\in\Delta_+$ and
$\frac{kNE_{kN\delta}}{q^{2kN}-q^{-2kN}}$,
$\frac{kNF_{kN\delta}}{q^{2kN}-q^{-2kN}}$ lead in the
limit $\hbar\to 0$ to derivations
of $\hsle$, which we denote by $e_\gamma$, $f_\gamma$, $\hat{e}_k$,
$\hat{f}_k$ correspondingly. The action of automorphism $\omega$  on these
derivations inherits from its action on corresponding root vectors.

\subsection{The Universal $R$-matrix at roots of 1}

Now let's consider the expression of universal $R$-matrix \re{rmatrix}
in the limit $\hbar\to 0$. The singularities, which appear in all
$q$-exponents, are the same type as in the expression of universal
$R$-matrix of $\slh$. New type singularities appear due to
the factor $\frac{kN}{q^{2kN}-q^{-2kN}}$ in the exponent before all term
$E_{kN\delta}\ot F_{kN\delta}$ for any natural $k$.

But as in $\sle$ case, the adjoint action $\hat{R}$ of $R$ on $\hsle\ot\hsle$
is well defined.

Indeed, the adjoint action of every $q$-exponent term
$$
R_\gamma=
\exp_\q((q-q^{-1})(E_\gamma\ot F_\gamma)), \quad
\gamma=\alpha_i+n\delta
$$
in \re{rmatrix} can be treated as it has been done in $\sle$ case:

\begin{eqnarray*}
\label{R}
\lim_{\hbar\to 0} \Ad(R_\gamma)&=&
	\prod_{m=0}^{N-1}\Ad((1-\ve^mE_\gamma\ot F_\gamma)^{-\frac{m}{N}})\\
	&\times& \exp\left(c_{N^\prime}
	\frac{\ln(1-E_\gamma^N\ot F_\gamma^N)}{E_\gamma^N\ot F_\gamma^N}
	(e_\gamma\ot F_\gamma^N+E_\gamma^N\ot f_\gamma)\right),
\end{eqnarray*}
where $c_{N^\prime}$ is defined by \re{cN}.

{}From \re{compar} and \re{Drin} it follows that the operations
$$
\hat{e}_k=\lim_{\hbar\to
0}\ad\left(\frac{kNE_{kN\delta}}{q^{kN}-q^{-kN}}\right), \
\hat{f}_k=\lim_{\hbar\to
0}\ad\left(\frac{kNF_{kN\delta}}{q^{kN}-q^{-kN}}\right)
$$
also are the derivations on $\hsle$, as it was mentioned above. So,
\begin{eqnarray*}
\hat{R}_{kN\delta}&=&\lim_{\hbar\to 0}\Ad(R_{kN\delta})=\lim_{\hbar\to 0}
	\Ad\left(\exp\left(\frac{kN}{q^{kN}-q^{-kN}}E_{kN\delta}\ot F_{kN\delta}
	\right)\right)\\
 &=&\exp({\hat{e}_k}\ot F_{kN\delta}+E_{kN\delta}\ot \hat{f}_k),
\end{eqnarray*}
gives rise to an outer automorphism of $\hsle$.

Finally, the right term in \re{rmatrix} has the following adjoint action
$$
\hat{{\cal K}}= \Ad(\ve^{\half{1}H_0\ot H_0+c\ot d+d\ot c}) =
	K_0^{1\ot\half{1}h_0}
	\ot K_0^{\half{1}h_0\ot 1}\ve^{c\ot\ad(d)+\ad(d)\ot c}
$$
Here $h_0=\ad(H_0)$ is a derivation on $\hslq$.

So, we proved, that the quantum algebra $\hsle$ is autoquasitriangular
Hopf algebra with the automorpism
\be
\label{aut}
\hat{R}=\left(\prod_{\gamma\in\Delta_+}\hat{R}_{\gamma}\right)\hat{K},
\ee
where the product over positive roots is ordered according to the normal
order \re{normal}.

\subsection{The universal $R$-matrix on Verma modules}

Consider now Verma module $M_{\hat{\lambda}}$ over $\hsle$ with
highest weight $\hat{\lambda}$. It is generated by vectors
$$
v_{k_1\dots k_n}^{\hat{\lambda}}=F_{\gamma_n}^{k_n}\ldots F_{\gamma_1}^{k_1}
	v_0^{\hat{\lambda}}
	\quad k_1,\dots,k_n=0,1,\dots
	\quad\gamma\in\Delta_+
	\quad \gamma_1<\dots<\gamma_n,
$$
where $v_0^{\hat{\lambda}}$ is a highest weight vector:
$$
E_\gamma v_0^{\hat{\lambda}}=0 \quad Hv_0^{\hat{\lambda}}=
	\hat{\lambda}(H)v_0^{\hat{\lambda}}
$$
As for $\slq$-case all terms $R_\gamma$ and $K$ in the product of
universal $R$-matrix \re{rmatrix} are well defined in the limit
$\hbar\to0$. Indeed, there is a well defined action of derivations
$e_i$, $\hat{e}_i$  on $M_{\hat{\lambda}}$  by
$$
e_igv_0^{\hat{\lambda}}:=e_i(g)v_0^{\hat{\lambda}} \qquad
\hat{e}_igv_0^{\hat{\lambda}}:=\hat{e}_i(g)v_0^{\hat{\lambda}} \
\quad \forall g\in\hsle
$$
Moreover, in the action of \re{rmatrix}, on any vector
$x\in M_{\hat\lo}\ot M_{\hat\lo}$  the term $R_\gamma$ with
sufficiently large $\gamma$ give rise to identity and only finite number
of $R_\gamma$ survive. In the decomposition of each such $R_\gamma$
also survive only finitely many terms. So, the action of $R$ on
$x\in M_{\hat\lo}\ot M_{\hat\lo}$ is well  defined

\vspace{.5cm}

To define the action of the  Universal $R$-matrix  \re{rmatrix} on
$\slq$-Verma modules, the spectral parameter dependent homomorphism $\rho_x$:
$\hslq \rightarrow \slq$ must be introduced \cite{J2}:

\begin{eqnarray*}
\label{}
\rho_x(E_{\alpha_0})=E & \rho_x(F_{\alpha_0})=F &  \rho_x(H_0)=H \\
\rho_x(E_{\alpha_1})=xF & \rho_x(F_{\alpha_1})=x^{-1}E &  \rho_x(H_1)=-H
\end{eqnarray*}

Note, that in this representation the central charge $c$ is zero.
Under the action of $\rho_x$ the root vectors acquire the form (\cite{Zh}):

\begin{eqnarray}
\label{hsltosl}
&E_{\alpha_0+n\delta}
	=(-1)^nx^nq^{-nh}E \qquad F_{\alpha_0+n\delta}=(-1)^nx^{-n}Fq^{nh}
\nonumber\\
&E_{\alpha_1+n\delta}
	=(-1)^nx^{n+1}Fq^{-nh} \qquad F_{\alpha_1+n\delta}=(-1)^nx^{-n-1}q^{nh}E
\nonumber\\
&E^\prime_{n\delta}=
	\frac{(-1)^{n-1}}{{[}2{]}_q}x^nq^{-(n-1)h}(EF-q^{-2}FE)
\label{eq:En}\\
&F^\prime_{n\delta}=\frac{(-1)^{n-1}}{{[}2{]}_q}x^{-n}q^{(n-1)h}(FE-q^{-2}EF)
\nonumber
\end{eqnarray}

Substituting this in the expression of affine universal $R$-matrix following
\cite{Zh}, one can obtain the spectral parameter $R$-matrix:
\be
\label{decompos}
R\left(\frac{x}{y}\right)=(\rho_x\ot\rho_y)R=
	R^{+}\left(\frac{x}{y}\right)
	R^{0}\left(\frac{x}{y}\right)R^-\left(\frac{x}{y}\right){\cal K},
\ee
where
\begin{eqnarray}
\label{spec}
R^+(z) &=&\prod_{n\ge 0}\exp_\q\left((q-q^{-1})z^n
	(q^{-nH}E\ot Fq^{nH})\right)
\nonumber\\
R^0(z) &=& \exp\left(\sum_{n>0}\frac{n}{q^{2n}-q^{-2n}}z^n
        E_{n\delta}\ot F_{n\delta}\right)
\\
R^-(z) &=& \prod_{n\ge 0}\exp_\q\left((q-q^{-1})z^{n+1}
	(Fq^{-nH}\ot q^{nH}E)\right) \nonumber\\
{\cal K} &=& q^{{1\over 2}H\ot H}\nonumber
\end{eqnarray}

Now we consider \re{spec} on Verma modules $M_\lambda$ of $\slq$ and its
behavior at roots of unity.

Note that one can represent the terms $R^\pm,R^0$ of universal
$R$-matrix in a more suitable way by performing infinite sum
and infinite product in \re{spec}. So, we have (\cite{Khor}):

\begin{eqnarray}
\label{R+}
R^+(z) &=&1+(E\ot F)\frac{(q-q^{-1})}{1-zq^{-2} K^{-1}\ot K}\nonumber\\
&+&\frac{(E\ot F)^2}{(2)_{q^{-2}}!}\frac{(q-q^{-1})^2}
	{(1-zq^{-2} K^{-1}\ot K)(1-zq^{-4} K^{-1}\ot K)}
\\
&\ldots& +
\frac{(E\ot F)^n}{(n)_{q^{-2}}!}\frac{(q-q^{-1})^n}
	{(1-zq^{-2} K^{-1}\ot K)\dots (1-zq^{-2n} K^{-1}\ot K)}
	+\ldots \nonumber
\end{eqnarray}

\begin{eqnarray}
\label{R-}
&&R^-(z) =1+\frac{z(q-q^{-1})}{1-zq^{-2} K^{-1}\ot K}F\ot E\nonumber\\
&&+\frac{1}{(2)_{q^{-2}}!}\frac{z^2(q-q^{-1})^2}
	{(1-zq^{-2} K^{-1}\ot K)(1-zq^{-4} K^{-1}\ot K)}(F\ot E)^2+\ldots
\\
&&+\frac{1}{(n)_{q^{-2}}!}\frac{z^n(q-q^{-1})^n}
	{(1-zq^{-2} K^{-1}\ot K)\dots (1-zq^{-2n} K^{-1}\ot K)}(F\ot E)^n
	+\ldots \nonumber
\end{eqnarray}
and
\be
\label{R0}
R^0(z)=f(z )\bar{R}^0(z)
\ee
 where
\begin{eqnarray}
\label{fz}
f(z )&=&\exp  \sum_{n\geq 1}
\left( (q-q^{-1})
\frac{[\lambda_1 n]_q[\lambda_2 n]_q}{[2n]_q}
	\right)\frac{z^n}{n}\nonumber\\
	&=&\frac{(zq^{\lambda_1-\lambda_2-2};q^{-4})_\infty
	(zq^{\lambda_1-\lambda_2-2};q^{-4})_\infty}
	{(zq^{\lambda_1+\lambda_2-2};q^{-4})_\infty
	(zq^{-\lambda_1-\lambda_2-2};q^{-4})_\infty} ,
\end{eqnarray}
$$
(z;q)_\infty=\prod_{i=0}^{\infty}(1-zq^k)
$$
\begin{eqnarray}
\label{R'}
&&\bar{R}^0(z)= \exp  \sum_{n\geq 1}\left(
\frac{q^n+q^{-n}}{(q^n-q^{-n})} \left( q^{-\lambda_1 n}- K^{-n}\right)
\otimes \left( K^n-q^{\lambda_2 n}\right)\right)
\frac{z^n}{n}
\\
&&\times\exp  \sum_{n\geq 1}
\left( (q^{-\l n}-K^{-n})\otimes
q^{-n}\frac{[\lambda_2 n]_q}{[n]_q}+
q^n\frac{[\lambda_1 n]_q}{[n]_q}
\otimes (K^n-q^{\lambda_2 n})\right)
\frac{z^n}{n}
\nonumber
\end{eqnarray}

By performing the infinite sum in \re{R'} one can easy show
that the  term $\bar{R}^0(z)$ acting on $v_i^\lo\ot v_j^\lt$
gives rise to the following expression, which is well defined in the
limit  $q^N\to 1$:

\begin{eqnarray}
\label{R'R''}
 \bar{R}^0(z)v_i^\lo \ot v_j^\lt &=&
	\frac{\prod_{l=j-i+1}^j(1-q^{-2l}q^{\lambda_2-\lambda_1}z)}
	{\prod_{l=j-i+1}^i(1-q^{2l}q^{\lambda_2-\lambda_1}z)}
\nonumber\\
& \times & \frac{\prod_{l=0}^{j-1}(1-q^{-2l}q^{\lambda_2+\lambda_1}z)}
	{\prod_{l=0}^{i-1}(1-q^{2l}q^{-\lambda_2-\lambda_1}z)}
	v_i^\lo \ot v_j^\lt
\end{eqnarray}

The scalar factor $f(z)$ \re{fz}  is singular for $q^{N^\prime}=1$. It can
be omitted from the expression of $R$-matrix.  So, the regular expression
of $R$-matrix for $q^N=1$ on $M_\lo\ot M_\lt$ has the form
\be
\label{Rverma}
R_{\lo,\lt}(z) = R^+(z)\bar{R}^0(z)R^-(z)
\ee
Note that it satisfy $R_{\lo,\lt}(z)v_0^\lo\ot v_0^\lt=v_0^\lo\ot
 v_0^\lt$. This renormalized expression of $R$-matrix does't satisfy
the quasitriangularity condition \re{quasi}. The intertwining
property \re{inter} and spectral parameter dependent Yang-Baxter
equation
\be
	R_{\lo,\lt}(\frac{x_1}{x_2})R_{\lo,\lambda_3}(\frac{x_1}{x_3})
	R_{\lt,\lambda_3}(\frac{x_2}{x_3})
	=R_{\lt,\lambda_3}(\frac{x_2}{x_3})R_{\lo,\lambda_3}(\frac{x_1}{x_3})
	R_{\lo,\lt}(\frac{x_1}{x_2})
\ee
are satisfied.

Let us consider now the possibility to restrict \re{Rverma} on finite
dimensional semicyclic modules. Recall that the semicyclic module
$V_{\alpha,\lambda}$ is obtained by factorisation of $M_\lambda$ on
$I_{\alpha,\lambda}=(F^N-\alpha) M_\lambda$
for some $\alpha\in \ab{C}$:
$$
V_{\alpha,\lambda}= M_\lambda / I_{\alpha,\lambda}
$$
The $R$-matrix \re{Rverma} is well defined on
$V_{\alpha_1,\lambda_1}\ot V_{\alpha_2,\lambda_2}$ if it preserves this
factorization, i.e.
\be
\label{factor1}
R_{\lo,\lt}(z) (M_{\lambda_1}\ot I_{\alpha_2,\lambda_2})\subset
	 (M_{\lambda_1}\ot I_{\alpha_2,\lambda_2})\bigoplus
	 (I_{\alpha_1,\lambda_1}\ot M_{\lambda_2})
\ee
and
\be
\label{factor2}
R_{\lo,\lt}(z) (I_{\alpha_1,\lambda_1}\ot M_{\lambda_2})\subset
	 (M_{\lambda_1}\ot I_{\alpha_2,\lambda_2})\bigoplus
	 (I_{\alpha_1,\lambda_1}\ot M_{\lambda_2})
\ee
The conditions above follow from
\begin{eqnarray*}
&& R_{\lo,\lt}\left({x\over y}\right) (\lt^N\cdot F^N\ot 1 +1\ot F^N)
\nonumber\\
&&\qquad =(\lo^N \cdot 1\ot F^N+F^N\ot 1) R_{\lo,\lt}\left({x\over y}\right)
\end{eqnarray*}
\begin{eqnarray*}
&&R_{\lo,\lt}\left({x\over y}\right)
	({x^N}\cdot F^N\ot 1 +y^N\lo^N\cdot 1 \ot F^N)\nonumber\\
&&\qquad	=(y^N\cdot 1\ot F^N+x^N\lt^N\cdot F^N\ot 1)
	R_{\lo,\lt}\left({x\over y}\right)
\end{eqnarray*}
Here we used the intertwining property \re{inter} for
\begin{eqnarray*}
\Delta(E_i^N)=E_i^N\ot 1+K_1^{-N}\ot E_1 & \Delta(F_i^N)=F_i^N\ot K_i^N+
	1\ot F_i^N
\end{eqnarray*}
So, one can express the operators
$$
R_{\lo,\lt}\left({x\over y}\right)(F^N\ot 1) \ \  {\rm and} \ \
R_{\lo,\lt}\left({x\over y}\right)(1\ot F^N)
$$
as a linear combination of
the operators
$$
(F^N\ot 1)R_{\lo,\lt}\left({x\over y}\right) \ \  {\rm and} \  \
(1\ot F^N)R_{\lo,\lt}\left({x\over y}\right)
$$
(if $\frac{x^N}{y^N}\ne\lo^N\lt^N$). In the same way,
$$
R_{\lo,\lt}\left({x\over y}\right)((F^N-\lo)\ot 1)  \ \  {\rm and} \ \
R_{\lo,\lt}\left({x\over y}\right)(1\ot (F^N-\lt))
$$
are a linear conbinations
of terms
$$
((F^N-\lo)\ot 1)R_{\lo,\lt}\left({x\over y}\right) \ \  {\rm and} \ \
(1\ot (F^N-\lt))R_{\lo,\lt}\left({x\over y}\right)
$$
with the same coefficients if parameters $x, y, \lo,
\lt, \alpha_1, \alpha_2$ lie on the algebraic curve
\begin{eqnarray}
\label{curve}
\frac{\alpha_1}{1-\lo^{N}}=\frac{\alpha_2}{1-\lt^{N}} &
	z^N=\left(\frac{x}y\right)^N=1
\end{eqnarray}
In this case the factorisation conditions \re{factor1},\re{factor2} are
fulfilled and $R$-matrix \re{Rverma} can be reduced to $R$-matrix
$R_{V_{\alpha_1,\lo}\ot V_{\alpha_2,\lt}}$ of semicyclic representations of
$\hsle$, considered in \cite{GS,IU,HS1}. The condition \re{curve}  on
parameters of
representations appears naturally as a consistency of factorisation
$V_{\alpha,\lambda}=M_\lambda / I_{\alpha,\lambda}$ with the intertwining
property \re{inter} of $R$-matrix.

Note that the formulae \re{R+},\re{R-},\re{R'},\re{Rverma} can be applied
directly to
semicyclic modules, using the constraint $F^N=\alpha\cdot {\rm id}$ on
$V_{\alpha,\lambda}$.

\section{Discussions}
Let's consider now the possibility of restriction of the automorphism \re{aut}
in evaluation representation \re{} to cyclic modules. Recall that their
intertwining operators  are Boltsman weight of Chiral Potts model (\cite{BS}).
The cyclic modules are representations of quotient algebra
$Q_\xi=Q_{\beta,\alpha,\lambda}$, $\xi=(\beta,\alpha,\lambda)$,
which is
obtained from $\sle$ by factorisation on ideal
$I_{\beta,\alpha,\lambda}$, generated by
$(F^N-\alpha), (E^N-\beta),
(K^N-\lambda^N)$, $(\beta,\alpha,\lambda\in\ab{C})$ (\cite{Kac}):
$$
Q_{\beta,\alpha,\lambda}= \sle / I_{\beta,\alpha,\lambda}
$$
The necessary condition for restriction of $\hat{R}(z)$ to
$Q_\xi$ is the constraint on parameters of representation
to lie on the algebraic curve, defined by
\begin{eqnarray}
\label{curve1}
&\frac{\alpha_1}{1-\lo^{N}}=\frac{\alpha_2}{1-\lt^{N}} \qquad
	\left(\frac{x}y\right)^N=1 & \nonumber\\
&\frac{\beta_1}{1-\lo^{-N}}=\frac{\beta_2}{1-\lt^{-N}} &
\end{eqnarray}

We expect that this condition is sufficient also and automorphism $\hat{R}$
can be restricted on some automorphism (outer, in general) of quotient algebra
$Q_{\xi_1}\ot Q_{\xi_2}$, which we denote by
$\hat{R}^{Q_{\xi_1}\ot Q_{\xi_2}}$.

Consider now its action on tensor product
of cyclic modules $V_{\xi_1}\ot V_{\xi_2}$.
$\hat{R}^{Q_{\xi_1}\ot Q_{\xi_2}}$ reduced here to matrix
algebra automorphism. Recall that every automorphism of matrix algebra is
inner. So,
$$
 \hat{R}^{Q_{\xi_1}\ot Q_{\xi_2}}_{V_{\xi_1}\ot V_{\xi_2}}=
	\Ad(R_{\xi_1,\xi_2})
$$
with some matrix $R_{\xi_1,\xi_2}$. This $R$-matrix is nothing but  to
the Boltsmann
weights of Chiral Potts model.

For quotients $Q_{0,\alpha,\lambda}$, corresponding
to semicyclic irreps, this suggestion is true.

Note that in case of $q^4=1$ there is a Hopf algebra homomorphism between
different quotients, as it was observed in \cite{A}.
This fact was used there to construct $R$-matrices of
quotient algebras for $q^4=1$ from the $R$-matrix of $Q_{0,0,\lambda_1}\ot
Q_{0,0,\lambda_1}$, which corresponds to nilpotent irreps.

Another question is to extend these results in case of other quantum
algebras.

Then we had finished this work, we saw the papers \cite{Pet} and
\cite{Kac1}, where
the center of quantum Kac-Moody algebras  was studied also. As it was
observed there the
automorphisms $\omega_\pm$ \re{sigma} corresponds to translations of
quantum Weyl group.

\vspace{1cm}

The authors thanks V.Pasquer, J-B.Zuber for helpfull discussions and
Saclay Theory Division for hospitality, where this work was started.

\end{document}